\documentstyle[12pt]{amsart}
\makeatletter
\def\LaTeX{\leavevmode L\raise.42ex
\hbox{\kern-.3em\size{\sf@size}{0pt}\selectfont A}\kern-.15em\TeX}
\makeatother

\newcommand{\BibTeX}{{\rm B\kern-.05em{\sci\kern-.025emb}\kern-.08em\TeX}}

\makeatletter
\def\@currentlabel{2.1}\label{e:dispaa}
\def\@currentlabel{2.21}\label{e:dispau}
\def\@currentlabel{2.22}\label{e:dispav}
\def\@currentlabel{2.23}\label{e:dispaw}
\def\@currentlabel{2.24}\label{e:dispax}
\def\theequation{\thesection.\@arabic\c@equation}
\makeatother

\newcommand{\M}{{\cal M}}

\newcommand{\RR}{ I\!\!R}

\newcommand{\N}{{ I\!\!N }}

\newcommand{\beq}{\begin{equation} }
\newcommand{\eqq}{\end{equation} }
\newcommand{\cuad}{{\sqcap\kern-.68em\sqcup}}

\newcommand{\equ}[1]{(\ref{#1})}

\renewcommand{\theequation}{\thesection.\arabic{equation}}

\newtheorem{definition}{Definition}[section]
\newtheorem{teo}{Theorem}[section]

\newtheorem{lema}{Lemma}[section]

\newtheorem{corollary}{Corollary}[section]
\newtheorem{remark}{Remark}[section]
\newcommand{\bremark}{\begin{remark} \em}
\newcommand{\eremark}{\end{remark} }

\newcommand{\ve}{\varepsilon}

\begin{document}
\title[Super-linear elliptic equation]{Super-linear elliptic equation for the Pucci operator without growth restrictions for the data}

\author{Maria J. Esteban}
\address{M.J. Esteban, Ceremade UMR CNRS 7534, Universit\'e Paris Dauphine, 75775 Paris Cedex 16, France}
\author{Patricio L. Felmer}
\address{\noindent P. Felmer, Departamento de
Ingenier\'{\i}a  Matem\'atica\\ and
Centro de Modelamiento Matem\'atico, UMR2071 CNRS-UChile\\
 Universidad de Chile, Casilla 170 Correo 3, Santiago, Chile.}
\author{Alexander Quaas}
\address{A. Quaas, Departamento de  Matem\'atica, Universidad T\'ecnica Federico Santa Mar\'{i}a
Casilla: V-110, Avda. Espa\~na 1680, Valpara\'{\i}so, Chile.
}


\begin{abstract}
In this paper we deal with existence and uniqueness of solution to super-linear problems for the Pucci operator:
$$\,-\M^+(D^2u)+|u|^{s-1}u=f(x)\,\quad \mbox{in } \RR^n,
$$ where $s>1$ and $f$ satisfies 
only local integrability conditions. This result is well known  when, instead of the Pucci 
operator, the Laplacian or a divergence form operator is considered. 
Our existence results use the   Alexandroff-Bakelman-Pucci inequality 
since we cannot use any variational formulation. For radially symmetric $f$ we
 can prove our results under less  local integrability assumptions, taking advantage of an appropriate variational formulation.
We also obtain  an existence result with boundary  explosion in smooth domains.
 
\end{abstract}
\date{}\maketitle

\setcounter{equation}{0}

\section{Introduction}

For parameters $0<\lambda\le\Lambda$ we consider $\M_{\lambda,\Lambda}^+$ 
and $\M_{\lambda,\Lambda}^-$, the maximal Pucci 
operators as defined in \cite{caffarelli2}. 
Whenever no confusion arises we will simply write $\M^+$ and $\M^-$, omitting 
the parameters. The problem we study in this article is the solvability of the 
differential equation 
\begin{equation}
-\M^+(D^2u)+|u|^{s-1}u=f(x)\qquad\mbox{in}\quad \RR^N,\label{main}
\end{equation}
when $s>1$ and for $f$ having only local properties, 
but without assuming any growth condition at infinity. 

When $\M^+$ is replaced by the Laplace operator, Brezis showed in \cite{brezis}
 that whenever $s>1$, one can find a (unique) solution to the above problem 
assuming only local integrability of $f$.  This very weak assumption is enough 
when the nonlinearity is increasing and super-linear, as in the case of $\, 
|u|^{s-1}u\,$ with $s>1$.
This result was extended to the case of a general quasilinear operator, including
the $p$-Laplace operator, and to parabolic equations by   
Boccardo, Gallouet and V\'azquez 
in \cite{boccardo1} and  \cite{boccardo3}, respectively. See also the work by Leoni
in \cite{Le} where more general nonlinearities are considered.
In all these works, the existence of the solution is obtained using in a crucial way
the variational structure of the equation by choosing appropriate test 
functions to obtain {\sl a priori } estimates. 

The Pucci operator is fully nonlinear and has no variational structure. 
So, in order to find a solution to \eqref{main}, we have to work in the 
viscosity solution framework and we cannot use test functions and integration 
by parts to derive {\sl a priori } estimates. 
The use of the viscosity theory forces us to work in the $L^N(\RR^N)$ 
framework and indeed, the presence of the $\, |u|^{s-1}u\,$ term in the 
equation allows us also to prove the existence of a unique $\,L^N$-viscosity 
for \eqref{main} whenever $\,f\in L^N_{loc}(\RR^N)$. Since there is no available
theory for viscosity solution when $f\in L^1_{loc}(R^N)$, at this point we cannot 
expect 
to obtain results under this weaker condition. 
However, in view of our results in Section \S 3 for the radially symmetric case, 
one may expect to find solutions
when
$f$ has less than $L^N$-integrability, but at this point we are not able to do it.
Our first theorem is the following
\begin{teo}
\label{teo1}
Assume that $s>1$.
For  every function $f\in L^N_{loc}(\RR^N)$,
the equation \equ{main}
possesses a unique solution in the $L^N$-viscosity sense and if $f\ge 0$ a.e. then
$u(x)\ge 0$ for all $x\in \RR^N$. A similar result holds if we replace $\M^+$ by 
$\M^-$ in \equ{main}. 
\end{teo}
The formal definition of solution is given in Section \S 2. 

\smallskip

It is well known that in the case of  super-linear problems one can find 
solutions which explode at the boundary of a bounded domain. This has 
been shown for various cases of linear and nonlinear second order elliptic 
operators in divergence form. See for instance the work by 
Keller \cite{keller}, Loewner and Nirenberg \cite{loewner}, Kondrat'ev and 
V. Nikishkin \cite{otros}, Diaz and Letelier \cite{diaz}, Diaz and Diaz \cite{diaz-diaz}, Del Pino and Letelier 
\cite{delpino} and Marcus and Veron \cite{marcus}.

In the case of the Pucci operator, the techniques used to prove Theorem \ref{teo1} can also be used to prove the 
following theorem on the existence of solutions in a bounded set, 
with explosion on the boundary. The simplest situation is the following
\begin{teo}
\label{teo3}
Assume that $s>1$,  $f\in L^N(\Omega)$ 
and let $\Omega\subset\RR^N$ be a bounded open set of class $C^2$.
Then the equation
\begin{eqnarray}
-\M^+(D^2u)+|u|^{s-1}u&=&f\qquad\mbox{in}\quad \Omega,\label{mainexplo1}\\
\lim_{x\to\partial\Omega}u(x)=\infty\label{inf}
\end{eqnarray}
possesses a solution in the $L^N$-viscosity sense.
\end{teo}
Here we only address the simplest situation, but the same kind of results should also hold true under more general assumptions. Moreover, the asymptotic study of the blow-up rate, both when $f \in L^N(\Omega)$ as above and when $f$ itself explodes at the boundary, is an interesting problem, due to the nonlinearity of the differential operator.

\medskip

In the second part of this paper we analyze the case of radially symmetric 
data $f$. Here  we can prove existence and uniqueness of solutions under
 weaker integrability assumptions on $f$.
The reason for this is that in the radial case we can
re-write equation \equ{main} as a divergence form quasilinear ordinary
differential equations, for which one can define a notion of weak solution. 
In this case we are back to integration 
by parts techniques.

Comparison between radial solutions and positivity results however, 
are not  obtained in a direct way.
This is because the coefficient 
of the second order derivative in the equation depends on the solution and its 
first derivative in a nonlinear way. Thus,  when comparing two solutions 
we do not have an obvious common factor for the second derivative of the 
difference or, if we have it we do not 
control its integrability at the origin.
An {\sl ad hoc} argument has to be found to do comparison in this case, see Lemma 3.3.
\begin{teo}\label{teo2I} 
Assume $s>1$ and $f$ is a radially symmetric function satisfying
\begin{equation}\label{condf}
\int_0^Rr^{N_+-1}|f(r)|dr<\infty,
\end{equation}
 for all $R>0$.
Here  $\,N_+:= \frac{\lambda}{\Lambda}(N-1)+1$, with $\,\lambda$ and $\Lambda\,$ being the parameters defining the Pucci operator  $\,\M_{\lambda,\Lambda}^+$.
Then equation \equ{main} has a unique weak radially symmetric solution 
and if  $f$ is nonnegative then $u$
 is also nonnegative.
\end{teo}

The formal definition of radially symmetric weak solution and the proof of Theorem 1.3
are given in Section \S 3. See also Remark 3.1 
 where we discuss  the assumptions on $f$ in this case.

\begin{remark}
In all our results, the power function $\, |u|^{s-1}u\,$ could be  replaced with 
nonlinear functions which are super-linear at infinity, however for simplicity 
all throughout the paper we will only deal with the pure power case.
In this direction see \cite{brezis}, \cite{boccardo1} and \cite{Le}.
Let us also stress that the assumption $s>1$ is essential for our results to
 hold, as we can see from the discussion in \cite{boccardo1}.\end{remark}

\setcounter{equation}{0}

\section{ The general case with $f\in L^N_{loc}(\RR^N)$.}

We devote this section to prove Theorem \ref{teo1} by an approximation procedure
 together with a local estimate based on a truncation argument and the application of 
the Alexandroff-Bakelman-Pucci inequality.

We start recalling the notion of solution suitable when  the right hand 
side in \equ{main} is only in $L^p_{loc}(\RR^N)$. Following the work by
Caffarelli, Crandall, Kocan and Swiech \cite{caffarelli1}, we first notice 
that the framework requires $p>N-\ve_0$, where $\ve_0>0$ depends on the ellipticity
constants $\lambda$ and $\Lambda$. Thus the case $p=N$, which is our framework is covered by the theory.
 Even though the context of the definitions
in \cite{caffarelli1} is much more general, for the purposes of this article  we  only consider a 'semilinear' case
\equ{main}
\begin{equation}
-\M(D^2u)+F(u)=f(x) \qquad\mbox{in}\quad \RR^N,\label{main1}
\end{equation}
where $\M$ stands for $\M^+$ or $\M^-$ and $F$ is an increasing continuous odd 
function.
Following 
\cite{caffarelli1}
we have the following definition: 
\begin{definition}
{\it Assume that $f\in L^p_{loc}(\RR^N)$, then  we say that a continuous function 
$u:\RR^N\to\RR$ is
an $L^p$-viscosity subsolution (supersolution) of the equation \equ{main1} in $\RR^N$
if for all $\varphi\in W^{2,p}_{loc}(\RR^N)$ and a point $\hat x\in \RR^N$ at which
$u-\varphi$ has a local maximum (respectively, minimum) one has
\begin{eqnarray}
  {\rm ess}\liminf_{x\to\hat x} (-\M(D^2\varphi(x))+F(u(x))-f(x) )\le 0  
\label{des1}
\end{eqnarray}
 \begin{eqnarray}     ({\rm ess}\limsup_{x\to\hat x} (-\M(D^2\varphi(x))+F(u(x))-f(x) )\ge 0 ).
 \label{des2}
\end{eqnarray}
Moreover, $u$ is an $L^p$-viscosity solution of \equ{main1} if it is both and $L^p$-viscosity subsolution and an $L^p$-viscosity supersolution.}
\end{definition}
\medskip

In what follows we say that $u$ is a $C$-viscosity (sub or super) solution of 
\equ{main1} when in the definition above we replace the tests function space
$\varphi\in W^{2,p}_{loc}(\RR^N)$ by $C^2(\RR^N)$. In this case the limits \equ{des1} and \equ{des2} become  simple
evaluation at $\hat x$, as given in \cite{crandall1}.

As we mentioned above, the idea is to consider a sequence of approximate problems and 
then take the limit at the end. So, given $f\in L^N_{loc}(\RR^N)$ we assume 
$\{f_n\}$ is a sequence of $C^\infty(\RR^N)$ functions so that for every bounded set $\Omega$
\begin{equation}\label{lim}
\lim_{n\to\infty}\int_\Omega|f_n-f|^Ndx =0.
\end{equation}
The sequence $\{f_n\}$ is easily constructed by mollification and a 
diagonal argument. 

The following is a basic existence and regularity result we need in our construction of a solution to \equ{main}.
\begin{lema}\label{lema1}
For every $n\in\N$ there is a solution $u_n\in C^2(B_n)$ of the equation
\begin{eqnarray}
-\M(D^2u_n)+\frac{1}{n} u_n +|u_n|^{s-1}u_n &=& f_n(x)\quad\mbox{in}\quad B_n \label{aprox1}
\\ 
u_n&=&0\quad\mbox{on}\quad \partial B_n,\label{aprox2}
\end{eqnarray}
where $B_n=B(0,n)$ is the ball centered at $0$ and with radius $n$. Here $\M$ stands for $\M^+$ or $M^-$.
\end{lema}

\noindent
{\bf Proof.}
We observe that there is a $M_n$ so that
$$
-M_n^s\le f_n(x)\le M_n^s\quad\mbox{for all}\quad x\in B_n
$$
and then $v_-=-M_n$ and $v_+=M_n$ are subsolution and subsolution of 
\equ{aprox1}-\equ{aprox2}, respectively. 
Then we can use the existence Theorem 4.1 in \cite{crandall1} for viscosity solutions of \equ{aprox1}-\equ{aprox2} to find $u_n$ a $C$-viscosity solution.
We observe that the hypothesis of Theorem 4.1 are fully satisfied by our operator, 
which is proper and satisfies the other hypothesis with $\gamma=1/n$, see 
\cite{crandall1}.

Noticing that $u_n$ solves the equation
\begin{equation}\label{En}
\M(D^2u_n)=g_n
\end{equation}
for the continuous function $g_n(x)=u_n(x)/n+|u_n(x)|^su_n(x)-f_n(x)$ we find that
$u_n\in C^{0,\alpha}(B_n)$ for $\alpha>0$, 
applying Proposition 4.10 in \cite{caffarelli2}.
Then we observe that $g_n$ is in $C^{0,\beta}(B_n)$, for certain $\beta>0$ and we may apply the regularity theory of Caffarelli \cite{caffarelli3} 
to obtain $u_n\in C^{2,\beta}(B_n)$.$\Box$

\medskip

Our next lemma is a version of Kato's inequality for $C$-viscosity
solutions  of equation
\equ{main1} with continuous right hand side.
\begin{lema}\label{lema2} 
Assume $\Omega\subset \RR^N$ and $u,v,f:\Omega \to\RR$ are continuous functions
and let $G(x)=F(u(x))-F(v(x))$.
If $u-v$ is a $C$-viscosity
solution  of equation
\begin{equation}
-\M(D^2 (u-v))+G(x)\le f \qquad\mbox{in}\quad \Omega\label{main5}
\end{equation}
then  $(u-v)^+$ is a $C$-viscosity  solution
of 
\begin{equation}
-\M(D^2 (u-v)^+)+G^+\le f^+ \qquad\mbox{in}\quad \Omega.\label{main6}
\end{equation}
Here $\M$ stands for $\M^+$ or $\M^-$.
\end{lema}
\noindent
{\bf Proof.}
If $x\in \Omega$ satisfies $u(x)-v(x)>0$ or $u(x)-v(x)<0$ then obviously $u-v$ 
satisfies \equ{main6} at $x$. If $u(x)-v(x)=0$ then we choose a test function 
$\varphi$ 
so that $(u-v)^+-\varphi$ has a local maximum at $x$, but then 
$(u-v)-\varphi$ has a local maximum at $x$ and then we may use \equ{main5} to obtain
$$
-\M(D^2\varphi(x))\le f^+
$$ 
so that \equ{main6} is satisfies in $x$, since $G(x)=0$.
$\Box$

Now we give a generalization of Kato's inequality (see \cite{Kato}) for $C$-viscosity
solutions  of equation
\equ{main1}.

\begin{lema}\label{remarku}
If we assume that  $u,f:\Omega \to\RR$ are continuous functions and $u$ is a
 $C$-viscosity
solution  of equation
\begin{equation}\label{YY1}
-\M(D^2 u)+F(u)= f \qquad\mbox{in}\quad \Omega,\label{main9}
\end{equation}
then  $|u|$ satisfies
\begin{equation}\label{YY2}
-\M^+(D^2 |u|)+F(|u|)\le |f| \qquad\mbox{in}\quad \Omega\label{main7}
\end{equation}
 in the $C$-viscosity sense.
\end{lema}
\noindent
{\bf Proof.}
In case $\M=\M^+$, we first use $v=0$ in Lemma \ref{lema2} to get that $u^+$ 
is a subsolution with $f^+$ as right hand side, and then observe that
$$
-\M^+(D^2(-u))+F(-u)\le f^-,
$$
since $\M^-\le \M^+$, that gives that $u^-$ is a subsolution with $f^-$ as right hand side. We conclude that
$|u|=\max\{u^+,u^-\}$ satisfies \equ{main7}

The case $\M=\M^-$ is similar.
$\Box$

\medskip

The following lemma contains the crucial local estimate 
for solutions of
\equ{aprox1}-\equ{aprox2} of class $C^2$. This result was proved 
by Brezis \cite{brezis} in the context of the Laplacian and tells that solutions 
have local estimates independent of the global behavior of $f$. The approach in
 \cite{brezis}, see also \cite{boccardo1}, is to use suitable test functions and 
integration by parts. This cannot be done here since the differential operator does not have
divergence form. For this 
result the fact that $s>1$ is essential.

\begin{lema}\label{lema3}
Let $s>1$ and $g$ continuous in $\Omega\subset \RR^N$, an open set. 
Suppose that $g\ge 0$ in $\Omega$ and $u$ is a 
nonnegative $C$-viscosity solution of 
$$
-\M(D^2u)+\frac{1}{n}u+|u|^{s-1}u\le g\qquad\mbox{in}\quad\Omega,
$$
 which additionally satisfies the  inequality in 
the classical sense whenever $u(x)>0$, then for all $R>0$ and $R'>R$ such that
$B_{R'}\subset\Omega$
\begin{equation}
\sup_{B_R}u\le C(1+\|g\|_{L^N(B_{R'})}),
\label{desigualdad}
\end{equation}
where $C=C(s,R,R',N,\lambda,\Lambda)$ does not depend on $g$ nor $n$.
Here $\M$ stands for $\M^+$ and $\M^-$.
\end{lema}

\noindent
{\bf Proof.}
We assume that $u$ is non trivial, otherwise the estimate is obvious.
In what follows we write $\Omega^+=\{x\in\Omega\,/\, u(x)>0\}$.
Let $\xi(x)=(R')^2-|x|^2$ and $\beta=2/(s-1)$ and consider $v=\xi^\beta u$. 

Consider the contact set for the function $v$, which is defined as 
$$\Gamma_v^+=\{x\in B_{R'}\,/\, \exists\, p\in\RR^N \mbox{ with  } 
v(y)\le v(x)+\langle p,y-x\rangle,\,\, \forall y\in B_{R'}\}.
$$
We observe that 
$\Gamma_v^+\subset \Omega^+\cap B_{R'}$ and that if $\bar v$ is the concave envelope
of $v$ in $\overline{B}_{R'}$ 
then for $x\in B_{R'}$ we have $v(x)=\bar v(x)$ if and only if 
$x\in \Gamma_v^+$. The function $\bar v$, being concave, satisfies
$$
\bar v(y)\le v(x)+\langle Dv(x),y-x\rangle,
$$
for all $x\in \Gamma_v^+$ and $y\in \overline{B}_{R'}$. Choosing adequately $y\in\partial
B_{R'}$ we obtain
\begin{equation}
|Dv(x)|\le\frac{v(x)}{R'-|x|}\qquad\mbox{for all}\quad x\in \Gamma_v^+\label{cota}
\end{equation}
and consequently,
\begin{equation}
\xi^\beta|Du(x)|\le\frac{v(x)}{R'-|x|}+\beta|D\xi|\xi^{-1} v(x) \qquad\mbox{for all}\quad x\in \Gamma_v^+.\label{cota0}
\end{equation}
Now we see that for all $x\in \Gamma_v^+$ the function $v$ satisfies
\begin{equation}\label{cota1}
-\M^+(D^2v)+\frac{1}{n}v+\xi^{\beta(1-p)}v^s\le \xi^\beta f+I+II+III,
\end{equation}
where 
\begin{equation}
I=-\beta\xi^{\beta-1}u\M^-(D^2\xi)\le C\xi^{-2}v,\label{cota2}
\end{equation}
\begin{equation}
II=-\beta(\beta-1)\xi^{\beta-2}u\M^-(D\xi\otimes D\xi)\le 
C\xi^{-2}v\label{cota3}
\end{equation}
and \begin{eqnarray}
III&=&-\beta\xi^{\beta-1}\M^-(D\xi\otimes Du+Du\otimes D\xi)\nonumber\\
&\le & c\xi^{\beta-1}|Du|\nonumber\\
&\le & c(R'+|x|+\beta |D\xi|)\xi^{-2}v\le C\xi^{-2}v,\label{cota4}\end{eqnarray}
where we used \equ{cota0}. Here $c$ and $C$ are constants depending on $R'$ and $s$.
From \equ{cota1}-\equ{cota4} we find that $v$ satisfies
$$
-\M^+(D^2v)+\xi^{-2}v(|v|^{s-1}-C)\le \xi^\beta g  
\qquad\mbox{for all}\quad x\in \Gamma_v^+.
$$
Now we define $w=\max\{v-C^{1/(p-1)},0\}$ in $B_{R'}$ and we observe that 
$\Gamma_w^+\subset \Gamma_v^+$ and 
$\Gamma_w^+\{x\in B_{R'}\,/\, w>0\}$. Consequently
$$
-\M^+(D^2w)\le \xi^\beta g,  \qquad\mbox{for all}\quad x\in \Gamma_w^+.
$$
Thus, from Alexandroff-Bakelman-Pucci inequality (see for example \cite{caffarelli2}),
$$
\sup_{B_{R'}}w\le C\|\xi^\beta g\|_{L^N(B_{R'})},
$$
but then
$$
\sup_{B_{R}}u\le \sup_{B_{R'}}v\le \sup_{B_{R'}}w +C^{1/(p-1)}\le
C(1+\|g\|_{L^N(B_{R'})}),
$$
where $C$ represents a generic constant depending only on $s,R,R',N,\lambda$ and
$\Lambda$ but not on $g$ nor $n$, as desired.

The case of $\M^-$ is similar. $\Box$
\begin{remark}
Observe that in this estimate the constant $C$ does not even depend on the possibly 
arbitrary values of $u$ on $\partial\Omega$. 
 This fact is very important in the
study solutions of this equation having explosion on the boundary of $\Omega$,
as we see in Theorem
\ref{teo3}. 
\end{remark} 
\noindent
{\bf Proof of Theorem \ref{teo1} (Existence)}
We start with a sequence of smooth functions $\{f_n\}$ such that for 
every bounded set 
$\Omega$ \equ{lim} holds. Then we use Lemma \ref{lema1} to 
construct a sequence of solutions $\{u_n\}$ of equation
\equ{aprox1}-\equ{aprox2}. According to Lemma \ref{remarku} and 
\ref{lema3}, for every $0<R<R'<n$ we have 
$$
\sup_{B_{R'}}|u_n|\le
C(1+\|f\|_{L^N(B_R)}),
$$
where $C$ does not depend on $f$ nor in $n$.
With this inequality in hand we look at equation \equ{En} and use Proposition 4.10 
in \cite{caffarelli2} to obtain, for every bounded open set $\Omega$,
$$
\|u_n\|_{C^\alpha(\Omega)}\le C,
$$
where $C$ does not depend on $n$, but only on  $f$, $\Omega$ and the other 
parameters.
By a diagonal procedure, we then obtain a subsequence
of
solutions  of equation
$$
-\M(D^2 u_n)+c_nu_n+|u_n|^{s-1}u_n=f_n, 
$$ 
 that we keep calling
$\{u_n\}$,  
such that $u_n$ converges uniformly
over every bounded subset of $\RR^N$. Here the equation holds in $B_{1/c_n}$,
with $c_n\to 0$ as $n\to\infty$, and $f_n$ has been renamed.
 Then using Theorem 3.8 in \cite{caffarelli1} we conclude that $u$ is an 
$L^N$-viscosity solution of \equ{main}, completing the proof of the existence part of the Theorem \ref{teo1}.
$\Box$

\medskip

The next lemma gives the positivity part of Theorem \ref{teo1}.
\begin{lema}\label{lemapos}
Assume $s>1$. If  $f\le 0$ a.e. and $u$ solves 
equation \equ{main} in the $L^N$-viscosity sense then $u\le 0$ in $\RR^N$ and
 if  $f\ge 0$ a.e. then  $u\ge 0$ in $\RR^N$.  
Similar results hold if we replace $\M^+$ by $\M^-$ in \equ{main}.\end{lema}
\noindent
{\bf Proof.}
We proceed as in \cite{brezis}, considering the function defined by Osserman in \cite{Osserman}:
$$
U(x)=\frac{CR^\alpha}{(R^2-|x|^2)^\alpha}\quad\mbox{in } B_R \quad R>0,
$$
where $\alpha=2/(s-1)$ and $C^{s-1}=2\alpha\Lambda\max\{N, \alpha+1\}$. Since $U'$ 
and $U''$ are positive, we see that $\M^+(D^2U)=\Lambda \Delta U$ and then a 
direct computation gives that
\begin{equation}\label{inn}
-\M^+(D^2U)+U^s\ge 0 \quad\mbox{in } B_R.
\end{equation}
From here, the equation for $u$ and the non-positivity of $f$ we obtain
$$
-\M^+(D^2(u-U))+|u|^{s-1}u-U^s\le 0.
$$
We observe that this inequality is in the $L^N$-viscosity sense, however
since $f$ was dropped, it also holds in the $C$-viscosity sense. Then by 
Lemma \ref{lema2} we find
$$
-\M^+(D^2(u-U)^+)+(|u|^{s-1}u-U^s)^+\le 0
$$
from where we get
$$
-\M^+(D^2(u-U)^+)\le 0 \quad\mbox{in } B_R.
$$
We observe that the function $u-U$ is negative in the set $R-\delta\le |x|<R$, for 
some sufficiently small $\delta>0$. Then by Alexandroff-Bakelman-Pucci maximum principle $(u-U)^+= 0$
which implies $u-U\le 0$ in $B_R$. From here, taking point-wise limit as 
$R\to \infty$ we find that $u\le 0$. 

In 
case $f \geq 0$ we proceed similarly, but relying in Lemma 2.2 with the operator
$\M^-$, to obtain  that $u+U\ge 0$ in $B_R$. From here the result follows.
The arguments for \equ{main} with $\M^-$ instead of $\M^+$ go along the same lines.
$\Box$

\medskip
\noindent
{\bf Proof of Theorem \ref{teo1} (Uniqueness).}
We only discuss the case of operator $\M^+$, since the other goes in a similar way.
If $u_1$ and $u_2$ are solutions of \equ{main} then the continuous function
$w=u_1-u_2$ satisfies
$$
-\M^+(D^2w)+|u_1|^{s-1}u_1-|u_2|^{s-1} u_2\le 0
$$
in the $C$-viscosity sense. If fact $u_1$ and $u_2$ are in $W^{2,N}$ and 
satisfy the equation in the a.e. and so  the inequality above,
and then in the $L^N$-viscosity sense, 
see Lemma 2.6 and Corollary 3.7 in \cite{caffarelli1}. Since the inequality
does not have $L^N$ ingredients it is satisfied in the $C$-viscosity sense.
Next we use Lemma \ref{lema2} to obtain $$
-\M^+(D^2w^+)+(|u_1|^{s-1}s u_1-|u_2|^{s-1} u_2)^+\le 0
$$
and using that
$$
||a|^{s-1}a-|b|^{s-1}b|\ge \delta|a-b|^s,\quad \forall a,b\in\RR,
$$ 
for certain $\delta>0$ we conclude that
\begin{equation}\label{inw}
-\M^+(D^2w^+)+\delta(w^+)^s\le 0.
\end{equation}
Using Lemma \ref{lemapos} we obtain that $u_1-u_2\le0$. Interchanging the roles of 
$u_1$ and $u_2$ we complete the proof.
$\Box$

\medskip

Next we give an existence theorem for explosive solutions, whose proof follows
easily from the estimate given in Lemma \ref{lema3} and the known results for
the Laplacian. We keep in the simplest form, but
we think it may be extended to more general situations in particular considering a 
blowing up right hand side.

\noindent
{\bf Proof of Theorem \ref{teo3}}. We first consider an increasing sequence of smooth functions
$\{f_n\}$ such 
$$
\lim_{n\to\infty}\int_\Omega|f_n-f|^N =0.
$$
Then we find $u_n$ a solution to the problem
\begin{eqnarray}
-\M^+(D^2u_n)+\frac{1}{n}u_n+|u_n|^{s-1}u_n&=&f_n\qquad\mbox{in}\quad \Omega,\label{mainexplo1n}\\
u_n&=&n\qquad\mbox{in}\quad \partial\Omega
\end{eqnarray}
By comparison theorem we obtain $u_{n+1}\ge u_n$ 
in $\Omega$ for all $n\in\N$. 
By arguments similar to those given in the proof of Theorem \ref{teo1} (Existence), using Lemma \ref{lema3},
we obtain a subsequence, we keep calling $\{u_n\}$ so that
$u_n$ converges uniformly to a solution $u$ of \equ{mainexplo1}. 
 Moreover
$ u\ge u_n$ in $\Omega$ for all 
$n$ so that $\liminf_{x\to \partial\Omega}u\ge n$, for all $n$, so that
$u$ also satisfies \equ{inf}.
$\Box$

\setcounter{equation}{0}

\section{ The radial case }

When we are dealing with radially symmetric functions then the Pucci 
operator has a much simpler form. 
Since the eigenvalues of $D^2u$ are $u''$ of multiplicity one and 
$u'/r$ with multiplicity $N-1$ and defining
$
\theta(s)=\Lambda$ if $s\ge 0$ and $\theta(s)=\lambda$ if $s< 0$,  
then we easily see that for every $u$ radially symmetric,
$$
\M^+(D^2u)(r)=\theta(u''(r))u''(r)+\theta(u'(r))(N-1)\frac{u'(r)}{r}.
$$
Then we see that equation \equ{main} in the classical sense becomes
\begin{equation}\label{mainradial}
-\theta(u''(r))u''-\theta(u'(r))(N-1)\frac{u'}{r}+|u|^{s-1}u=f(r),
\end{equation}
for a radial function $f$.
In order to write this equation in a more simple form,
 we make some definitions. 
First we observe that for solutions of \equ{mainradial} we have
$$
\theta(u''(r))=\theta\{-\theta(u'(r))(N-1)\frac{u'}{r}+|u|^{s-1}u-f(r)\},
$$ 
which is more convenient as we see. We define 
$$
\Theta(r,u(r),u'(r))=\theta\{-\theta(u'(r))(N-1)\frac{u'}{r}+|u|^{s-1}u-f(r))\},
$$
the 'dimension'
$$
N(r,u(r),u'(r))= \frac{\theta(u'(r))}{\Theta(r,u(r),u'(r))}(N-1)+1
$$
and the weights
$$
\rho(r,u(r),u'(r))=e^{\int_{1}^r\frac{N(\tau,u(\tau),u'(\tau))-1}{\tau}d\tau}
$$
and 
$$
\tilde\rho(r,u(r),u'(r))=\frac{\rho(r,u(r),u'(r))}{\Theta(r,u(r),u'(r))}.
$$
If we define 
$$
N_+=\frac{\lambda}{\Lambda}(N-1)+1,
\quad\mbox{and}\quad N_-=\frac{\Lambda}{\lambda}(N-1)+1
$$
we see that $N_+\le N(r,u(r),u'(r))\le N_-$ and also,
$$
r^{N_--1}\le \rho(r,u(r),u'(r))\le r^{N_+-1} \,\,\mbox{if}\,\, 0\le r\le 1
$$
and
$$
\frac{\rho}{\Lambda}\le \tilde \rho \le \frac{\rho}{\lambda}.
$$
With these definitions we find that \equ{mainradial} is equivalent to
\begin{equation}
-(\rho u')'+\tilde\rho |u|^{s-1}u=\tilde\rho f(r)\label{main1radial}.
\end{equation}
When no confusion arises we omit the arguments in the functions $\rho$ and $\tilde
\rho$, in particular when we write $\rho v'$ we mean 
$\rho(r,v(r),v'(r))v'(r)$ and so on.
What is interesting about equation \equ{main1radial} is that it
allows to define a weaker notion of solution which extends the $L^N$-viscosity
sense to more general $f$. 
With this new notion we can prove a theorem for the existence of radial solutions of 
\equ{main} with a weaker condition on $f$ than in the non-radial
case of Section \S 2. See Remark 
\ref{comparacion}.

We consider the set of test functions defined as 
$$
H=\{\varphi:[0,\infty)\to\RR\,\,/\,\, \exists\, \phi\in W^{1,\infty}_0(\RR^N) \mbox{
such that }\phi(x)=\varphi(|x|)
\},$$
where $W^{1,\infty}_0(\RR^N)$ denotes the space of functions in 
$ W^{1,\infty}(\RR^N)$ with compact support. 
\begin{definition}
We say that $u:[0,R]\to\RR$ is a weak solution of 
\equ{main1radial} with Dirichlet boundary condition at $r=R$, 
if $u$ is absolutely 
continuous in $(0,R]$, $u(R)=0,$
\begin{equation}\label{cond1}
\int_0^R\rho|u|^sdr<\infty,\quad \int_0^R\rho|u'|dr<\infty
\end{equation}
and
\begin{equation}\label{cond2}
\int_0^R\rho u'\varphi'+\tilde\rho|u|^{s-1}u\varphi dr=\int_0^R\tilde\rho f\varphi dr\quad 
\forall\varphi\in H.
\end{equation}
\end{definition}
Now we state our theorem precisely which is a more complete version of Theorem \ref{teo2I}
\begin{teo}\label{teo2} 
Assume $s>1$ and $f$ is a radial function satisfying for all $R>0$,
\begin{equation}\label{condf1}
\int_0^Rr^{N_+-1}|f(r)|dr<\infty.
\end{equation}
Then equation \equ{main1radial} has a unique weak solution $u$ and if  $f$ is nonnegative then $u$
 is also nonnegative.

Additionally, for any $1<q<2s/(s+1)$
\begin{equation}\label{cero2}
\int_0^r\rho|u'|^qdr<\infty \quad \mbox{for all}\quad R>0.
\end{equation}
Moreover, the function $\rho u'$ is differentiable a.e. in $(0,\infty)$ and consequently satisfies 
\begin{equation}\label{cero}
\lim_{r\to 0}(\rho u')(r)=0,\quad \lim_{r\to 0}\int_0^r\rho|u'|dr=0.
\end{equation}
\end{teo}
In order to prove the theorem above we will perform an approximation procedure as in the general case. 
Because the problem is radial and has a divergence form  formulation we can 
get better estimates and pass to the limit, under weaker assumptions on $f$.

By regularizing $f$ and using a diagonal procedure we may find a sequence of radial smooth functions
$\{f_n\}$
 such that for all $0<R$
\begin{equation}\label{condfn}
\lim_{n\to\infty}\int_0^Rr^{N_+-1}|f_n(r)-f(r)|dr=0.
\end{equation}
Moreover, we may assume that there exists a function $g:(0,\infty)\to\RR$ such that
$|f_n(r)|\le g(r)$ for all $r>0$ and $\int_0^Rr^{N_+-1}|g(r)|dr<+\infty$, for all $R>0$.

First we have an existence result for the approximate problems.
\begin{lema}\label{lema21}
For every $n$ there is solution $u_n$ in $C^2[0,n]$ satisfying $u_n(n)=0$, 
\equ{cond1} with $R=n$
and
\begin{equation}\label{cond2n}
\int_0^n\rho_n u_n'\varphi'+\tilde\rho_n
(c_nu_n+
|u_n|^{s-1}u_n)\varphi=
\int_0^n\rho_n f_n\varphi,\quad \forall\varphi\in H.
\end{equation}
where $\rho_n(r)=\rho(r,u_n(r),u_n'(r))$ (similarly for $\tilde\rho_n$).
 and
$\{c_n\}$ is a positive sequence converging to zero.
\end{lema}
\noindent
{\bf Proof.} We may use the same argument of Lemma \ref{lema1}
together with Da Lio and Sirakov 
symmetry result \cite{sirakov}.
$\Box$ 

\medskip

Now we get some estimates following the ideas of Boccardo, Gallouet and Vazquez in
\cite{boccardo1}.
\begin{lema}\label{lema22}
Let $\{u_n\}$ be the sequence of solutions found in Lemma \ref{lema21}. Then, for all $0<R$ and $m\in (0,s-1)$ 
there is a constant $C$
depending on $R,m,s,N,\lambda$ and $\Lambda$, but not on $f$ nor $n$, such that
for all $n\in\N$ we have
\begin{equation}\label{cotarho}
\int_0^R\rho_n|u_n|^sds \le C(1+
\int_0^{2R}r^{N_+-1}|f|dr)
\end{equation}
and 
\begin{equation}\label{cotarho1}
 \int_0^{2R}\frac{\rho_n |u_n'|^2dr}{(1+|u_n|)^{m+1}}\le C(1+
\int_0^Rr^{N_+-1}|f|dr).
\end{equation}
\end{lema}
\noindent
{\bf Proof.} We consider the function $\phi$ defined as
$$
\phi(t)=\int_0^t\frac{dt}{(1+s)^{m+1}},\,\, t\ge 0
$$
and extended as an odd function to negative $t$, which is smooth and bounded.
We also consider a cut-off function $\theta:[0,\infty)\to\RR$ being smooth, with 
support in $[0,2R]$, equal $1$ in $[0,R]$, $0\le \theta\le 1$ and $|\theta'|\le 2/R$.

We define $v=\phi(u)\theta^\alpha$, where $\alpha>2s/(s-1-m)$. Omitting the index $n$ in what follows, using
 $v$ as a test function we obtain
\begin{eqnarray}
\nonumber & & \int_0^{2R}\frac{m\rho|u'|^2\theta^\alpha dr}{(1+|u|)^{1+m}}+
\int_0^{2R}\tilde\rho |u|^{s-1}u\phi(u)\theta^\alpha dr\\
 &\le &
\int_0^{2R}\tilde\rho f \phi(u)\theta^\alpha dr-
\alpha\int_0^{2R}\rho u' \phi(u)\theta^{\alpha-1} \theta' dr\label{eee}
\\ &\le & C(\int_0^{2R} r^{N_+-1})|f|dr+\int_0^{2R} \rho|u'|\theta^{\alpha-1}dr),
\label{ineq1}\end{eqnarray}
where we drop the term with $c_n$ in the first inequality. Using Young inequality, for some $\varepsilon>0$, we have
\begin{eqnarray}
\int_0^{2R}\rho|u'|\theta^{\alpha-1}dr&\le& \varepsilon\int_0^{2R}
\frac{m\rho|u'|^2\theta^\alpha}{(1+|u|)^{1+m}}dr \nonumber\\
&+&
\frac{1}{4\varepsilon}\int_0^{2R} \rho(1+|u|)^{1+m}\theta^{\alpha-2}dr
\label{ineq2}
\end{eqnarray}
and again
\begin{eqnarray}
\int_0^{2R} \rho(1+|u|)^{1+m}\theta^{\alpha-2}dr
&\le& 
\varepsilon^2 \int_0^{2R} \rho(1+|u|)^{s}\theta^{\alpha}dr
\nonumber\\
&+&
\frac{C}{\varepsilon^2}\int_0^{2R} \rho \theta^{\frac{\alpha(s-m-1)-2s}{s-m-1}}dr\nonumber\\
&\le&C(\varepsilon^{-2}+\varepsilon^2 \int_0^{2R} \rho|u|^{s}\theta^{\alpha}dr),
\label{ineq3}
\end{eqnarray}
where $C$ is a generic constant independent of $\ve$. Here we used our choice of $\alpha$.

Next we observe that $|t|^s\le |t|^{s-1}t\phi(t)/\phi(1) +1$ for all $t\in\RR$. 
Using this in \equ{ineq3} and then using what one gets and \equ{ineq2} in 
\equ{ineq1}, with the choice of a sufficiently small $\varepsilon$ we finally obtain
the desired inequalities. $\Box$

\medskip

\begin{corollary}
For all   $q\in(1,2s/(s+1))$ and for every $0<R$ there is a constant as 
in Lemma \ref{lema22}, 
such that
\begin{eqnarray}
\int_0^{R}{\rho_n |u_n'|^q}dr \le C(1+
\int_0^{2R}r^{N_+-1}|f|dr).\label{crucial}
\end{eqnarray}
\end{corollary}
\noindent
{\bf Proof.}
By H\"older inequality
we find
\begin{eqnarray}\nonumber
\int_0^R{\rho |u'|^q}dr\le 
\left(\int_0^R \frac{\rho |u'|^2 dr}{(1+|u|)^{1+m}}  \right)^{\frac{q}{2}}\left(\int_0^R 
\rho(1+|u|)^{\frac{q(1+m)}{2-q}}dr\right)^{\frac{2-q}{2}},
\end{eqnarray}
then by our choice of $m$
in Lemma \ref{lema22}   it is possible to choose $q>1$ such that 
$(m+1)q/(2-q)< s$ and then from Lemma \ref{lema22} we obtain the result. With the adequate choice of $m$ we can cover the range of $q$.
$\Box$ 

\medskip

\noindent
{\bf Proof of Theorem \ref{teo2} (Existence).}
We consider the sequence of $\{u_n\}$ of solution found in Lemma \ref{lema21} satisfying \equ{cond2n}. In what follows we show that this sequence  converges to a weak
solution 
of \equ{mainradial}. 

Now, considering the estimates
in Lemma \ref{lema22}, we see that the function $\rho_nu_n'$ has weak derivatives in
any interval of the form $(r_0,R_0)$ with $0<r_0<R_0$. Since the function
$\rho_n$ is differentiable a.e., we obtain then that $u_n$ is twice 
differentiable a.e. and $u_n''$ is in $L^1(r_0,R_0)$, because of the equation
satisfied by $u_n$ and estimates in Lemma \ref{lema22}. 
From here we conclude that $u_n'$ and $u_n$ 
are uniformly bounded in $(r_0,R_0)$. 
Using the equation again we conclude then that $u_n''$ is bounded by an $L^1$ function in $(r_0,R_0)$, which implies that $u_n'$ is equicontinuous. By the Arzela-Ascoli 
Theorem there exists a differentiable function $u$ in the interval $(r_0,R_0)$
such that, up to a subsequence,  
$u_n$ and $u_n'$ converges to $u$ and $u'$ respectively, in a uniform way in the
interval $(r_0,R_0)$.

We may repeat this argument for any interval $(r_0,R_0)$, so that by a diagonal
procedure, we can prove that up to a subsequence, $\{u_n\}$ and $\{u_n'\}$ 
converge point-wise to a differentiable function $u:(0,\infty)\to\RR$.
Notice that $\{\rho_n\}$ converges point-wise to $\rho(r)=\rho(r,u(r),u'(r))$.

Next we use the estimate \equ{crucial}, to prove that the sequence
$\{\rho_nu_n'\}$ is equi-integrable in $[0,R]$ and then it converges in $L^1[0,R]$ to 
$\rho u'$, for all $R>0$. It is only left to prove that $\{\tilde\rho_n|u_n|^s\}$ converges
in $L^1[0,R]$. 
For this purpose we introduce, as in \cite{boccardo1},  a new 
function $\phi$ in $\RR$ defined as $\phi(\nu)=
\min\{\nu-t,1\}$ if $\nu\ge 0$ and extended as an odd function to all $\RR$, for a parameter $t>0$.
Then we consider inequality \equ{eee} with the cut-off function $\phi(u_n)\theta$ 
to get
$$
\int_{E_n^{t+1}\cap (0,R)}\tilde\rho_n |u_n|^{s} dr
 \le 
\int_{E_n^t\cap (0,2R)} \tilde\rho_n| f_n|  dr+
C\int_{E_n^t\cap (0,2R)} \rho_n |u_n'|  dr,
$$
where $E_n^t=\{r>0\,/\, |u_n(r)|>t\}$. From \equ{cotarho} and \equ{crucial} it follows
that the second integral approaches zero if $t\to\infty$. From here the 
equi-integrability of $\rho_n |u_n|^{s}$ follows and we conclude.

Finally \equ{cero} is consequence of the integrability properties just 
proved for $u_n$ that also hold for $u$. This finishes
the proof.
$\Box$

\medskip

Now we prove the remaining part of Theorem \ref{teo2}, that is uniqueness and
non-negativity of weak solutions. For this purpose it would be natural to use comparison arguments, 
however those are a bit delicate in this case. In fact, in a natural way we 
may define the notion of weak subsolutions (supersolution) by writing $\le$ 
and use only nonnegative
test functions in \equ{cond2}. It happen that, if $u$ is a weak subsolution and $v$ is a 
weak supersolution, we cannot be sure that $w=u-v$ is a weak subsolution, since
we do not have good control of $\rho w'$ at the origin.

We first consider 
the no-negativity 
of solutions of when $f$ is nonnegative. For that purpose we need to find appropriate
test functions.
\begin{lema}\label{lemaposr}
If $u$ is a solution of \equ{mainradial} in the weak sense and $f\le 0$ a.e. in $[0,\infty)$
then $u\le 0$ for all $r>0$.
\end{lema}
\noindent
{\bf Proof.}
As in the general case, we 
 consider the function
$U$ given in the proof of Lemma \ref{lemapos}, which satisfies \equ{inn} in $B_R$. 
On the other hand by the
regularity of $u$ given above,  we have that $u(x)=u(r)$ satisfies equation
\equ{main}
a.e.  We may subtract the equations for $U$ and $u$ and get 
$$
-\M^+(D^2(u-U))+|u|^{s-1}u-U^s\le 0\quad a.e.\,\, \mbox{in } B_R.
$$
If we write $w=u-U$ then we see that
\begin{equation}\label{cond5}
-(\rho w')'+\tilde\rho (|u|^{s-1}u-U^s)  \le 0 \quad\mbox{in}\,\, (0,R) \, a.e.
\end{equation}
Here the function
$\rho$ and $\tilde\rho$ are defined in the natural way with $\theta(r)
=\theta(w'(r))$ and 
$\Theta$ given by
$$
\Theta(r)=\theta(w''(r)) \quad \mbox{in}\,\, (0,R) \, a.e.
$$

We see that the function $w$ is negative near $R$. If there exists $0<r_1<r_2<R$
such that $w>0$ in $(r_1,r_2)$ and $w(r_1)=w(r_2)=0$ then we may choose the
function $\varphi$, defined as 
$\varphi=w$ in $(r_1,r_2)$ and $\varphi\equiv 0$ elsewhere, as a test function in
\equ{cond5} to get 
$$
\int_{r_0}^R\rho |w'|^2 +\tilde\rho (|u|^{s-1}u-U^s)w dr\le 0.
$$
But each term in the left hand side is positive, then $w=0$ in $(r_1,r_2)$.

Thus, either $w(r)\le 0$ in $(0,R)$ 
or
there is $r_0\in (0,R)$ such that $w>0$ in $(0,r_0)$ and $w(r_0)=0$. 
To see that the second  case is impossible we just need to prove that 
\begin{equation}\label{cero1}
\int_0^{r_0}\rho(w)|w'|dr<\infty \quad\mbox{and}\quad\lim_{r\to 0}(\rho(w) w')(r)=0,
\end{equation}
since in this case we may use the function $\varphi$, defined as
$\varphi=w$ in $(\bar r, r_0)$ and $\varphi\equiv w(\bar r)$ in $(0, \bar r)$,
as a test function in \equ{cond5} and get a contradiction.

Assuming \equ{cero1} for the moment,  we see that $u\le U$ in $[0,R]$ and this is true
 for all $R>0$. Taking
limit as $R$ goes to infinity, keeping $r$ fixed, we conclude that $u\le 0$ in
$[0,\infty)$. 

To complete the proof we show \equ{cero1}.
To see this, we first observe that there is $\bar r\in (0,r_0)$ such that $w'(\bar r)<0$
and then from inequality \equ{cond5} we find that $w''(r)>0$ a.e and $w'(r)<0$ in
$r\in (0,\bar r)$. A posteriori we see that $w''(r)>0$ a.e and $w'(r)<0$ in 
$r\in (0, r_0)$ and consequently $\rho(w)=r^{N_+-1}$ there.
Next we assume that $u'$ is negative at some point in $(0,r_0)$, because otherwise
the functions $u$ and $u'$ would be  bounded and then $w$ and $w'$ are bounded, 
	yielding \equ{cero1}.
Since $u''>U''>0$ in $(0,r_0)$ we see then that $u'<0$ near the origin and 
consequently $\rho(u)=r^{N_+-1}$.
Since \equ{cero} holds we see that \equ{cero1} holds.
$\Box$

\medskip
\noindent
{\bf Proof of Theorem \ref{teo1} (Uniqueness).}
Let $u_1$ and $u_2$ be two solutions of equation \equ{mainradial} in the weak sense, then they satisfies \equ{main} a.e. in $\RR^N$, with abuse of notation $u_i(x)=u_i(|x|)$,
$i=1,2$.  
Then we define $w=u_1-u_2$ and proceed as in the Proof of Theorem \ref{teo1} 
to obtain that $w$ satisfies \equ{inw} a.e. in $\RR^N$. Now we follows the proof of Lemma \ref{lemaposr}.
$\Box$

\begin{remark}\label{comparacion}
Let us consider a continuous function $f$ in 
$\RR^N\setminus\{x_i\,/\,i=1,...,k\}$, such that near each singularity
$$
f(x)\sim \frac{c_i}{|x-x_i|^{\alpha_i}},\quad x\sim x_i,\,\,i=1,...,k.
$$ 
In order to apply Theorem \ref{teo1} we need $\alpha_i<1$ for all
$i=1,...,k$.
In contrast, assuming that $f$ is radially symmetric with a singularity at the 
origin of the form
$$
f(r)\sim \frac{c}{r^{\alpha}}, \quad r\sim 0, r>0,
$$
in order to apply Theorem \ref{teo2}, we only need $\alpha<N_+$. 
We observe that if $\lambda/\Lambda\to 0$ then $N_+\to 1$, while if 
$\lambda/\Lambda=1$ then $N_+=N$.  

When we have a radial function $f$ being in $L_{loc}^p(\RR^N)$ with
$p>N/N_+$ then $f$ satisfies our hypothesis \equ{condf1} and we may apply Theorem \ref{teo2}. This is particularly interesting if $N$ and $N_+$ are close to each other.
\end{remark}
\begin{remark}
Let $f$ be a function in $\RR^N$ and define
$$
g(r)=\max\{|f(x)|\,/\, |x|=r\}
$$
and assume that $g$ satisfies \equ{condf1}. This will be the case if $f$ has a singularity of the form $r^{-\alpha}$ with $\alpha<N_+$.

Then we may construct a solution of \equ{main1radial}. This solution is a 'candidate'
for a supersolution for equation \equ{main} with $f$ as a right hand side. 
However, since the two notion of solutions are not compatible, this is not posible.
\end{remark}
\begin{remark}
In this section we have considered only the case of the Pucci operator $\M^+$, 
however these results can be adapted for the operator $\M^-$ as well.
\end{remark}
\noindent
{\bf Acknowledgements:}
The second author was  partially supported by Fondecyt Grant \# 1070314 
 and FONDAP de Matem\'aticas Aplicadas and
the third author
was partially supported by Fondecyt Grant \# 1040794 and Proyecto Interno USM \# 12.05.24 .
This research was partially supported by Ecos-Conicyt project C05E09.


\end{document}